# Renewal Theory and Geometric Infinite Divisibility

**E Sandhya**   and   **R N Pillai**

Department of Statistics　　　　　　　　　　　　　　　　Valiavilakom
Prajyoti Niketan College　　　　　　　　　　　　　　　　　Ookode
Pudukkad, Trichur – 680 301　　　　　　　　　Trivandrum – 695 020
**India.**　　　　　　　　　　　　　　　　　　　　　　　　**India.**

*esandhya@hotmail.com*

*Abstract.* The role of geometrically infinitely divisible laws in renewal equations and super position of renewal processes are explored here.

**AMS Classification (2000)**: 60E07, 60E10, 60G55, 60G07.

## 1. Introduction

The notion of geometrically infinitely divisible (GID) laws is by Klebanov *et al.* (1984). The idea was studied by Sandhya (1991*a*) in her thesis, where-in she had further developed: (i) geometric version of an ID law (Pillai and Sandhya, 2001), the notions of geometric attraction and partial geometric attraction (Sandhya and Pillai, 1999) (ii) its relation to distributions with completely monotone derivative (Pillai and Sandhya, 1990), *p*-thinning of renewal processes (Sandhya, 1991*b*), reliability theory (Pillai and Sandhya, 1996) and a generalization of GID laws to $\nu$-ID laws (Sandhya, 1996). This paper also is part of the thesis, Sandhya (1991*a*). For recent reviews on this topic see Kozubowski and Rachev (1999) and Satheesh (2003).

Following distributions are relevant here. Characteristic functions (CF) $\varphi(t) = 1/[1+\psi(t)]$, where $\psi(o) = 0$, $\psi(t) = a\psi(bt)$ for some $0<b<1<a$, $ab^\alpha = 1$, $0<\alpha \leq 2$, defines semi-$\alpha$-Laplace laws of order *b*, Pillai (1985). This family, also known as geometrically semi-stable, is invariant under geometric summation, Sandhya (1991*a*). Restricting the support of these

Visit: http://ernakulam.sancharnet.in/probstat/

laws to [0,∞) we get semi-Mittag-Leffler Laws with Laplace transform (LT) $\phi(s) = 1/[1+\psi(s)]$, where $\psi(o) = 0$, $\psi(s) = a\psi(bs)$ for some $0<b<1<a$, $ab^\alpha = 1$, where $0<\alpha \leq 1$, Sandhya (1991*b*). Distributions with CF $1/(1+\lambda|t|^\alpha)$, $0<\alpha \leq 2$ are α-Laplace (Pillai, 1985), and those with LT $1/(1+s^\alpha)$, $0<\alpha \leq 1$ are Mittag-Leffler laws, Pillai (1990). They form subclasses of semi-α-Laplace and semi-Mittag-Leffler laws respectively.

Pillai (1985) in theorem.1, showed that a d.f $F$ satisfies the equation

$$F(x) = pF(x/b) + qF*F(x/b), \text{ for some } 0<p<1, \text{ and } p+q=1 \qquad (1)$$

iff $F$ is semi-α-Laplace of order $b = p^{1/\alpha}$. This equation is nothing but a renewal equation with terminating process as solution. Here we will consider certain variations of this equation in section.2.

The superposition of two point processes $N_1$ and $N_2$ is the point process $N$ defined by $N_t = N_{1,t} + N_{2,t}$, $t \geq 0$. Here we assume that no two events of $N_1$ and $N_2$ occur at the same time, that is, there are no multiplicities. An earliest result on superposition due to Mc Fadden (1962) is: If $N = N_1 + N_2$ is a stationary renewal process and $N_2$ is a Poisson process independent of $N_1$, and if $N_1$ is a stationary renewal process with finite variance for the inter arrival times, then $N_1$ also is a Poisson process. Cinlar and Agnew (1968) consider the superposition of non-independent processes. They use an indicator process, which labels the points of the superposition process according to the components to which they belong.

Using the ideas of indicator process and *p*-thinning we improve upon certain results on superposition of renewal processes in section.3. In section.4 we show that certain examples in Feller (1968), are GID.

## 2. Renewal Equations and Geometric Infinite Divisibility

We begin with some corollaries to theorem.1 of Pillai (1985).

**Corollary.2.1** A non-negative distribution $F$ satisfies the renewal equation (1) for some $0<p<1$, iff $F$ is semi-Mittag-Leffler with $b = p^{1/\alpha}$.

**Corollary.2.2** A d.f $F$ on **R** (on **R**$^+$) satisfies (1) for every $0<p<1$, iff it is α-Laplace (Mittag-Leffler).

*Proof.* Follows from theorem.1 of Pillai (1985) by noticing that when the equation is satisfied for every $0<p<1$, the semi-$\alpha$-Laplace (semi Mittag-Leffler) law becomes the corresponding $\alpha$-Laplace (Mittag-Leffler) law.

As a variation of this we have:

**Corollary.2.3** An $\alpha$-Laplace law is the solution to $Z$ in the equation

$$Z(x) = z(x) + Z*F(x) \qquad (2)$$

where $F(x) = qF_\alpha(x/p^{1/\alpha})$ and $z(x) = pF_\alpha(x/p^{1/\alpha})$ for every $0<p<1$, $F_\alpha$ being the d.f of an $\alpha$-Laplace law. Similarly;

**Corollary.2.4** A Mittag-Leffler law is the solution to $Z$ in the renewal equation

$$Z(x) = z(x) + Z*F(x) \qquad (3)$$

where $F(x) = qF_\alpha(x/p^{1/\alpha})$ and $z(x) = pF_\alpha(x/p^{1/\alpha})$ for every $0<p<1$, $F_\alpha$ being the d.f of a Mittag-Leffler law.

**Theorem.2.1** The solution $Z(x)$ of the renewal equation (3) with $z(x) = pG$ for any $p \in (0,1)$ and $F$ replaced by $qG$, $p+q=1$, where $G$ is a proper d.f, is GID iff $G$ is GID.

*Proof.* Let $\varphi$ and $\phi$ are the LTs of $G$ and $Z$ respectively. Then (3) is equivalent to:

$$\phi = p\varphi + q\varphi\phi. \qquad (4)$$

Suppose $G$ is GID. Then $\varphi = 1/[1+\psi]$, $\psi(o) = 0$ and $\psi(s)$, $s>0$ has completely monotone derivative. Hence;

$$\phi = p/[1+\psi] + q\phi/[1+\psi] \text{ and subsequently}$$

$$\phi = 1/[1+\psi/p].$$

Now since $\psi$ has completely monotone derivative $\phi$ is GID.

Conversely, suppose that the solution of (3) is GID. That is;

$$\phi = 1/[1+\psi_1], \psi_1(o) = 0 \text{ and } \psi_1(s), s>0 \text{ has completely monotone}$$
derivative. Now a substitution in (4) gives

$$\varphi = 1/[1+p\psi_1] \text{ and hence } \varphi = \text{ is GID. That completes the proof.}$$

**Remark.2.1** With respect to equation (2), one can now see that Theorem.2.1 is true in general for GID laws on **R** because if $1/[1+\psi]$ is a CF then $1/[1+c\psi]$ is again a CF for any $c>0$.

**Corollary.2.5** In (4) if $G$ is gamma with LT $1/[1+s]^{\alpha}$. Then the solution is GID iff $\alpha \leq 1$.

**Corollary.2.6** In (4) if $G$ is has LT $1/[1+s^{\alpha}]^{\beta}$. Then the solution is GID iff $\alpha \leq 1$ and $\beta \leq 1$.

See Sandhya (1991*b*) for an interpretation of the last two corollaries in the context of subordination of random processes.

## 3. Superposition of Renewal Processes and Geometric Infinite Divisibility

Let $N_1$, $N_2$ be two point processes and let $N = N_1 + N_2$. If $t_1$, $t_2$, ….are the points of $N$, we define

$X_n = k$  iff  $t_n \in N_k$, $k = 1,2$.

The process $X = \{X_1, X_2, ….\}$ is called the indicator process for $(N_1, N_2)$. Cinlar and Agnew (1968) proved; Suppose that $X$ is a Bernoulli process independent of $N$. Then $N_1$ and $N_2$ are independent only if $N$ is a Poisson process (possibly non-stationary).

Here we are observing the problem in a different manner. We have; $N$ a renewal process and $X(p)$ which is distributed as Bernoulli($p$) independent of $N$. Given $N$, the process $N_1$ ($N_2$) may be thought of as that process obtained by $p$-thinning ($q$-thinning) of $N$ w.r.t $X(p)$.

Now consider $N$ as a renewal process where the inter arrival times are distributed as Mittag-leffler with LT $\varphi(s) = 1/[1+s^{\alpha}]$, $0<\alpha \leq 1$. Then the LT of the $p$-thinned process is $\varphi_p(s) = 1/[1+s^{\alpha}/p]$ which is again Mittag-Leffler (with a scale change). Here since $\varphi(s)$ is GID it serves as the inter arrival time distribution of a Cox and renewal process (Sandhya, 1991*b*) and this implies that $\varphi_p(s)$ (and similarly $\varphi_q(s)$) also corresponds to a renewal process and they are also Mittag-Leffler. Hence;

**Theorem.3.1** Suppose that $X(p)$ is Bernoulli($p$), independent of $N$ where $N$ is renewal with Mittag-Leffler as the inter arrival time distribution. Then $N_1$ and $N_2$ are also Mittag-Leffler renewal processes (with a scale change).

**Theorem.3.2** Suppose that $N$ is a renewal process and $X(p)$ is the indicator process, which is Bernoulli($p$) independent of $N$. Then $N_1(p)$ ($N_2(p)$) is renewal w.r.t $X(p)$ for some $p \in (0,1)$ having inter arrival time distribution of the same type as that of $N$ iff $N$ has an inter arrival time distribution which is semi-Mittag-Leffler.

*Proof.* Suppose that $N$ has a GID inter arrival time with LT $\phi(s)$. Now $N_1(p)$ has an inter arrival time distribution with LT

$$\phi_1(s) = p\phi(s)/[1 - q\phi(s)] \text{ for some } p \in (0,1)$$

and similarly the LT of $N_2(p)$ is $\phi_2(s) = q\phi(s)/[1 - p\phi(s)]$ for some $q \in (0,1)$. Now $\phi(s)$ and $\phi_1(s)$ or $\phi(s)$ and $\phi_2(s)$ are of the same type iff $\phi(s)$ is semi-Mittag-Leffler which follows from corollary.2.1.

**Corollary.3.1** in the above setup; $N_1(p)$ ($N_2(p)$) is renewal w.r.t $X(p)$ for all $p \in (0,1)$ having inter arrival time distribution of the same type as that of $N$ iff $N$ has an inter arrival time distribution which is Mittag-Leffler.

*Proof.* Follows from corollary.2.2.

 Recall: Mecke (1968) (as recorded in Yannaros (1989)) showed that a point process is Cox iff it can be obtained by *p*-thinning, for every thinning parameter $p \in (0,1)$. Yannaros (1989) proved that if a renewal process is a Cox process, then it cannot be the thinning of some non-renewal process, and all the possible original processes are *p*-thinnings of other renewal processes for every thinning parameter $p \in (0,1)$.

**Theorem.3.3** Suppose that $N$ is a renewal process and $X(p)$ is the indicator process, which is Bernoulli($p$) independent of $N$. Then $N_1(p)$ ($N_2(p)$) is renewal w.r.t $X(p)$ for all $p \in (0,1)$ iff $N$ has an inter arrival time distribution that is GID.

*Proof.* Suppose that $N$ has a GID inter arrival time with LT $\phi(s)$. Now the LTs of the inter arrival time distributions of $N_1(p)$ and $N_2(p)$ are respectively

$\varphi_1(s) = p\varphi(s)/[1 - q\varphi(s)]$ for all $p \in (0,1)$ and

$\varphi_2(s) = q\varphi(s)/[1 - p\varphi(s)]$ for all $q \in (0,1)$.

Hence $\varphi_1(s)$ ($\varphi_2(s)$) is the LT of a GID law and thus the inter arrival time distribution of a Cox and renewal process.

Conversely, assume that $N_1$ and $N_2$ are renewal with respect to $X_p$ for all $p \in (0,1)$. Thus $\varphi_1(s)$ and $\varphi_2(s)$ are LTs of GID laws and hence $N_1$ and $N_2$ are Cox and renewal. Hence by Yannaros' result $N$ is also renewal and Cox. Hence $\varphi(s)$ is GID, thus completing the proof.

**Theorem.3.4** Suppose $N$ is a point process and $X(p)$ is the indicator process that is independent of $N$. If $N_1(p)$ is Cox and renewal, then $X(p)$ is Bernoulli($p$) and $N$ and $N_2(p)$ are also Cox and renewal.

*Proof.* Follows from Yannaros' and Mecke's result.

Again from Mecke's result we have:

**Theorem.3.5** Suppose $N$ is a point process and $X(p)$ is the indicator process that is independent of $N$. If $N_1$ is Cox, then $X(p)$ is Bernoulli($p$) and $N_2$ also is Cox.

**Theorem.3.6** Suppose $N$ is a point process and $X(p)$ is the indicator process that is independent of $N$ and $N(p) = N_1(p) + N_2(p)$ for all $p \in (0,1)$. Then $N_1$ and $N_2$ are both Cox.

**Remark.3.1** Notice that here we do not need the assumptions of stationarity, independence and finite variance on the component process and we have discussed cases other than Poisson.

## 4. Some Examples in renewal theory that are GID

Here we discuss certain examples in renewal theory where the LTs correspond to GID laws. The examples are taken from Feller (1968).

**Example.4.1** Let $u_n = P\{E$ occurs at the $n^{th}$ trial$\}$ and

$f_n = P\{$the event $E$ occurs for the first time at the $n^{th}$ trial$\}$.

Let the event $E$ be transient, that is; $\sum_{i=1}^{\infty} f_i < 1$.

If $\varphi(s)$ and $\phi(s)$ are the LTs respectively of $f_n$ and $u_n$, then we have:

$$\phi(s)[1- \varphi(s)] = 1 \text{ , see Feller (1968, p.311)}.$$

Setting $\varphi(1) = f = q$, $1- q = p$ and $\varphi(s)/q = \omega(s)$ we have:

$$p\phi(s) = 1/[1+\tfrac{q}{p}(1- \omega(s))] = 1/[1+\psi].$$

Now $\psi(o) = 0$ and $\psi(s)$ has completely monotone derivative. Hence $p\phi(s)$ is GID.

**Example.4.2** For the example.4(*b*) in Feller (1968, p.313) consider the LT

$$U(s) = 1/\sqrt{(1 - 4pqe^{-2s})}. \text{ Here } |p\text{-}q|\ U(s) \text{ is GID.}$$

**Example.4.3** For the example.4(*c*) in Feller (1968, p.315) consider the LT

$$\overline{U}(s) = \frac{1-\sqrt{(1 - 4pqe^{-2s})}}{2pqe^{-2s}}. \text{ Here } p\overline{U}(s) \text{ is GID.}$$

**Acknowledgement.** Authors wish to thank the referee for pointing out a mistake and the suggestions that has improved the paper. First author's work was supported by a research fellowship from the CSIR, India, during 1988-1991. The ideas were presented in 1989 at the National Symposium on Probability and Statistics, hosted by the Centre for Mathematical Sciences, Trivandrum, India.